\documentclass{article}
\usepackage{amsmath,amsthm,amscd,amssymb}
\usepackage[mathscr]{eucal}
\usepackage{amssymb}
\usepackage{latexsym}

\theoremstyle{definition}
\newtheorem{Def}{Definition}[section]
\newtheorem{Thm}[Def]{Theorem}
\newtheorem{Prop}[Def]{Proposition}
\newtheorem{Rem}[Def]{Remark}

\newtheorem{Cor}[Def]{Corollary}

\newtheorem{Lem}[Def]{Lemma}

\numberwithin{equation}{section}

\begin{document}

\title{Congruences for Hermitian modular forms of degree $2$}

\author{Toshiyuki Kikuta}
\maketitle




\begin{abstract}
We give two congruence properties of Hermitian modular forms of degree $2$ over $\mathbb{Q}(\sqrt{-1})$ and $\mathbb{Q}(\sqrt{-3})$. The one is a congruence criterion for Hermitian modular forms which is generalization of Sturm's theorem. Another is the well-definedness of the $p$-adic weight for Hermitian modular forms. 
\end{abstract}








\section{Introduction and Results}
\label{intro}
\subsection{Congruence criterion}
\label{subsec:1-1}
Sturm \cite{St} studied congruence properties of elliptic modular forms. In \cite{St}, the number of coefficients which are required to check the coincidence modulo a prime of two modular forms is determined. In \cite{Ch-Ch1}, Choi and Choie studied the analog of Sturm's theorem for Jacobi forms. Choi and Choie~\cite{Ch-Ch2}, Poor and Yuen~\cite{Po-Yu} and the author~\cite{Ki} have independently obtained some generalizations of Sturm's theorem in the case of Siegel modular forms of degree $2$. The first aim of this paper is to generalize Sturm's theorem to the case of Hermitian modular forms of degree $2$.

We state our results more precisely. Let $\boldsymbol{K}$ be the imaginary quadratic number field $\mathbb{Q}(\sqrt{-1})$ or $\mathbb{Q}(\sqrt{-3})$, ${\mathcal O}_{{\boldsymbol K}}$ the ring of integers in ${\boldsymbol K}$, $d_{\boldsymbol K}$ the discriminant of ${\boldsymbol K}$ and $U_2({\mathcal O}_{{\boldsymbol K}})$ the Hermitian modular group of degree $2$ defined as 
\begin{align*}
U_2({\mathcal O}_{{\boldsymbol K}}):=\left\{M \in M_{4}({\mathcal O}_{{\boldsymbol K}})\;\big{|}\; ^t\overline{M} J_2 M =J_2 \right\},\quad J_2:=\begin{pmatrix} O_2 & -1_2 \\ 1_2 & O_2 \end{pmatrix}.
\end{align*}
We take a character $\nu _k$ on $U_2({\mathcal O}_{{\boldsymbol K}})$ as 
\begin{align*}
{\nu_k=}
\begin{cases}
{\det}^{k/2} & \text{for}\;\; \boldsymbol{K}=\mathbb{Q}(\sqrt{-1}),\\
{\det}^k & \text{for}\;\; \boldsymbol{K}=\mathbb{Q}(\sqrt{-3}).
\end{cases}
\end{align*} 
For a subring $R\subset \mathbb{C}$, we denote by $M_k^{(s)}(U_2({\mathcal O}_{{\boldsymbol K}}),\nu _k)_{R}$ the space of symmetric Hermitian modular forms of weight $k$ with character $\nu _k$ whose all Fourier coefficients belong to $R$. We define $\Lambda _2({\boldsymbol K})$ as 
\begin{align*}
\Lambda _2({\boldsymbol K}):=\{H=(h_{ij})\in Her_2({\boldsymbol K})\;|\; h_{ii}\in \mathbb{Z},\;\sqrt{d_{{\boldsymbol K}}}h_{ij}\in {\mathcal O}_{\boldsymbol K}\; \}, 
\end{align*}
where $Her_2({\boldsymbol K})$ is the set of all Hermitian matrices of size $2$ whose all components are in ${\boldsymbol K}$. For $f\in M_k(U_2({\mathcal O}_{\boldsymbol K}),\nu _k)_{\mathbb{C}}$, then we write $f=\sum _{H}a_f(H)q^H$ the Fourier expansion of $f$, where $q^{H}:=e^{2\pi i{\rm tr}(HZ)}$, $Z$ is an element of the Hermitian upper-half space of degree $2$, $H$ runs over all elements of semi-positive definite of $\Lambda _2({\boldsymbol K})$. Let $p$ be a prime and $\mathbb{Z}_{(p)}$ the local ring of all $p$-integral rational numbers. We have the following theorem. 
\begin{Thm}
\label{ThmM}
Let $k$ be an even integer, $p$ a prime with $p\ge 5$. In the case ${\boldsymbol K}=\mathbb{Q}(\sqrt{-1})$ (resp. ${\boldsymbol K}=\mathbb{Q}(\sqrt{-3})$), if $f\in M_k^{(s)}(U_2({\mathcal O}_{\boldsymbol K}),\nu _k)_{\mathbb{Z}_{(p)}}$ satisfies that $a_f\left(\begin{pmatrix} m & \alpha  \\ \overline{\alpha } & n \end{pmatrix}\right)\equiv 0$ mod $p$ for all $m$, $n$, $\alpha $ with $0\le m\le (5k+8)/40$, $0\le n\le (15k+4)/120$, $mn-N(\alpha )\ge 0$ (resp. with $0\le m\le (5k+9)/45$, $0\le n\le (10k+3)/90$, $mn-N(\alpha )\ge 0$), then we have $f\equiv 0$ mod $p$. 
\end{Thm}

\begin{Rem}
In fact, to show the assumption of the theorem, we need only to check the case $n\le m$. 
\end{Rem}

By Theorem~{\ref{ThmM}}, we obtain immediately the followings. 
\begin{Cor}
Let $k$ be an even integer, $p$ a prime with $p\ge 5$. Let ${\boldsymbol K}=\mathbb{Q}(\sqrt{-1})$ (resp. ${\boldsymbol K}=\mathbb{Q}(\sqrt{-3})$). \\
(1) If $f\in M_k^{(s)}(U_2({\mathcal O}_{\boldsymbol K}),\nu _k)_{\mathbb{Q}}$ satisfies that $a_f\left(\begin{pmatrix} m & \alpha  \\ \overline{\alpha } & n \end{pmatrix}\right)\in \mathbb{Z}_{(p)}$ for all $m$, $n$, $\alpha $ with $0\le m\le (5k+8)/40$, $0\le n\le (15k+4)/120$, $mn-N(\alpha )\ge 0$ (resp. with $0\le m\le (5k+9)/45$, $0\le n\le (10k+3)/90$, $mn-N(\alpha )\ge 0$), then we have $f\in M_k(U_2({\mathcal O}_{\boldsymbol K}),\nu _k)_{\mathbb{Z}_{(p)}}$. \\
(2) If $f\in M_k^{(s)}(U_2({\mathcal O}_{\boldsymbol K}),\nu _k)_{\mathbb{Q}}$ satisfies that $a_f\left(\begin{pmatrix} m & \alpha  \\ \overline{\alpha } & n \end{pmatrix}\right)=0$ for all $m$, $n$, $\alpha $ with $0\le m\le (5k+8)/40$, $0\le n\le (15k+4)/120$, $mn-N(\alpha )\ge 0$ (resp. with $0\le m\le (5k+9)/45$, $0\le n\le (10k+3)/90$, $mn-N(\alpha )\ge 0$), then we have $f=0$. 
\end{Cor}

\subsection{Well-definedness of $p$-adic weight}
\label{subsec:1-2}

Serre \cite{Se} defined the notion of $p$-adic modular forms and applied it to the construction of $p$-adic $L$ functions. Some mathematicians have attempted to generalize the theory of Serre's $p$-adic modular forms to the case of several variables~\cite{Bo-Na,Ich,Na0,Mu-Na}. In particular, Ichikawa~\cite{Ich} showed that the $p$-adic weight is well-defined for Siegel modular forms of general degree. The second aim of this paper is to show the well-definedness of the $p$-adic weight for Hermitian modular forms of degree $2$ over $\mathbb{Q}(\sqrt{-1})$ and $\mathbb{Q}(\sqrt{-3})$.

In the situation as subsection~\ref{subsec:1-1}, we have the following theorem. 
\begin{Thm}
\label{ThmM3}
Let $k$ and $k'$ be even integers and $p$ a prime with $p\ge 5$. If $f\in M_k^{(s)}(U_2({\mathcal O}_{\boldsymbol K}),\nu _k)_{\mathbb{Z}_{(p)}}$ and $g\in M_{k'}^{(s)}(U_2({\mathcal O}_{\boldsymbol K}),\nu _{k'})_{\mathbb{Z}_{(p)}}$ satisfy that $f\not \equiv 0$ mod $p$ and $f\equiv g$ mod $p^l$, then we have $k\equiv k'$ mod $(p-1)p^{l-1}$. 
\end{Thm}
This theorem indicates that the $p$-adic weight of Hermitian modular forms is ``well-defined''. We shall explain what the ``well-defined'' means. Let $p$ be a prime, $v_p$ the normalized additive valuation on $\mathbb{Q}_p$ (i.e. $v_p(p)=1$) and $g=\sum _{0\le H\in \Lambda _2({\boldsymbol K})}b(H)q^H$ the formal power series such that $b(H)\in \mathbb{Q}_p$ for all $H$. We call $g$ a $p$-$adic$ $Hermitian$ $modular$ $forms$ $of$ $degree$ $2$ if there exists a sequence of Hermitian modular forms $\{f_m\}$ of weight $k_m$ with character $\nu _k$ (which are not necessarily symmetric) such that 
\begin{align*}
\lim _{m\to \infty }f_m =g\quad (p-{\rm adically}),
\end{align*} 
in other words, 
\begin{align*}
\inf\{v_p(a_{f_m}(H)-b(H))\:|\:0\le H\in \Lambda _2({\boldsymbol K})\}\to \infty \ (m\to \infty ). 
\end{align*}
For simplicity, we consider only the case $p\ge 5$. We define a group ${\boldsymbol X}$ as ${\boldsymbol X}:=\varprojlim \mathbb{Z}/(p-1)p^{l-1}\mathbb{Z}$ $(=\mathbb{Z}_p\times \mathbb{Z}/(p-1)\mathbb{Z}$). We find that $\mathbb{Z}\subset {\boldsymbol X}$ by the natural imbedding. The following property follows from our theorem. 
\begin{Cor}
Let $\{f_m\in M_{k_m}^{(s)}(U_2({\mathcal O}_{\boldsymbol K}),\nu _{k_m})_{\mathbb{Q}}\}$ and $\{f'_m\in M_{k'_m}^{(s)}(U_2({\mathcal O}_{\boldsymbol K}),\nu _{k'_m})_{\mathbb{Q}}\}$ be two sequences of symmetric Hermitian modular forms of even weight. If both $\{f_m\}$ and $\{g_m\}$ converge $p$-adically to a same formal power series $g$, then there exist limits of the sequences $\{k_m\}$ and $\{k'_m\}$ in $\boldsymbol X$ and we have 
\begin{align*}
\lim _{m\to \infty }k_m=\lim _{m\to \infty }k'_m\ {\rm in}\ {\boldsymbol X}. 
\end{align*}
\end{Cor}
By this corollary, we can define a weight of the $p$-adic Hermitian modular form $g$ satisfying $g({}^tZ)=g(Z)$, similarly to the case of elliptic modular forms. In particular, we see that the $p$-adic limit of Hermitian Eisenstein series of weight $k_m$ depends only on the limiting value of $k_m$ in ${\boldsymbol X}$.


\section{Hermitian modular forms}
\label{sec:2}
\subsection{Definition and notation}
\label{sec:2.1}
The Hermitian upper half-space of degree $n$ is defined by
\[
\mathbb{H}_n:=\{ Z\in M_n(\mathbb{C})\;|\; \tfrac{1}{2i}(Z-{}^t\overline{Z})>0\;  \}
\]
where ${}^t\overline{Z}$ is the transposed complex conjugate of $Z$. The space
$\mathbb{H}_n$ contains the Siegel upper half-space of degree $n$
\[
\mathbb{S}_n:=\mathbb{H}_n\cap Sym_n(\mathbb{C}).
\]
Let $\boldsymbol{K}$ be an imaginary quadratic number field with discriminant $d_{\boldsymbol{K}}$ and ring of integers $\mathcal{O}_{\boldsymbol{K}}$. The Hermitian modular group $U_n({\mathcal O}_{{\boldsymbol K}})$ acts on $\mathbb{H}_n$ by the fractional transformation
\[
\mathbb{H}_n\ni Z\longmapsto M<Z>:=(AZ+B)(CZ+D)^{-1},
\;M= \begin{pmatrix} A & B \\ C & D \end{pmatrix}\in U_n(\mathcal{O}_{\boldsymbol{K}}).
\]
The subgroup $SU_n(\mathcal{O}_{\boldsymbol{K}}):=U_n(\mathcal{O}_{\boldsymbol{K}})\cap SL_{2n}(\boldsymbol{K})$ coincides with the full group $U_n(\mathcal{O}_{\boldsymbol{K}})$ unless $d_{\boldsymbol{K}}=-3$ or $-4$. Let $\nu$ be an Abelian character on $U_n(\mathcal{O}_{\boldsymbol{K}})$. We denote by $M_k(U_n(\mathcal{O}_{\boldsymbol{K}}),\nu)$ the space of Hermitian modular forms of weight $k$ and character $\nu$ for $U_n(\mathcal{O}_{\boldsymbol{K}})$. Namely, it consists of holomorphic functions $f:\mathbb{H}_n\rightarrow \mathbb{C}$ satisfying
\[
f|_k M(Z):={\det}(CZ+D)^{-k}f(M<Z>)=\nu (M)\cdot f(Z),
\]
for all $M= \begin{pmatrix} A & B \\ C & D \end{pmatrix}\in U_n({\mathcal O}_{\boldsymbol K}) $. The cusp form $f\in M_k(U_n(\mathcal{O}_{\boldsymbol{K}}),\nu)$ is characterized by the condition
\begin{align*}
\Phi (f|_k \begin{pmatrix}{}^t\overline{U} & 0 \\ 0 & U \end{pmatrix}) \equiv 0\quad \text{for}\;
\text{all}\;
U\in GL_n(\boldsymbol{K})
\end{align*}
where $\Phi$ is the Siegel $\Phi$-operator. A modular form $f\in M_k(U_n(\mathcal{O}_{\boldsymbol{K}}),\nu)$ is called $symmetric$ if $f({}^tZ)=f(Z)$. We denote by $M_k^{(s)}(U_n(\mathcal{O}_{\boldsymbol{K}}),\nu)$ the subspace consisting of symmetric modular forms. 

\subsection{Fourier expansion}
\label{sec:2.2}
If $f\in M_k(U_n({\mathcal O}_{\boldsymbol K}),\nu)$ satisfies the condition
\begin{align*}
f(Z+B)=f(Z) 
\quad \text{for}\;\text{all}\;
B\in Her_n(\mathcal{O}_{\boldsymbol{K}}),
\end{align*}
then $f$ has a Fourier expansion of the form
\begin{align*}
f(Z)=\sum_{0\leq H\in\Lambda_n(\boldsymbol{K})}a_f(H)e^{2\pi i\text{tr}(HZ)}
\end{align*}
where
\[
\Lambda_n(\boldsymbol{K}):=\{ H=(h_{ij})\in Her_n(\boldsymbol{K})| h_{ii}\in\mathbb{Z},\sqrt{d_{\boldsymbol{K}}}h_{ij}\in\mathcal{O}_{\boldsymbol{K}} \}.
\]
Put $\omega:=\frac{1}{2}(d_{\boldsymbol{K}}+\sqrt{d_{\boldsymbol{K}}})$ and define the 
matrices
$\dot{Z}=(\dot{z}_{ij})$ and $\ddot{Z}=(\ddot{z}_{ij})$ by
\[
\dot{Z}:=\frac{\omega{}^tZ-\bar{\omega}Z}{\omega-\bar{\omega}},\quad
\ddot{Z}:=\frac{Z-{}^tZ}{\omega-\bar{\omega}}.
\]
Then the above $f$ can be considered as a function of the $\frac{1}{2}n(n-1)$ complex
variables $\ddot{z}_{ij}$\,$(i<j)$ in $\ddot{Z}$ and of the $\frac{1}{2}n(n+1)$ complex variables $\dot{z}_{ij}$\,$(i\leq j)$ in $\dot{Z}$. Moreover, $f$ has period $1$ for each of these variables. If we define
\[
\dot{q}_{ij}:=\text{exp}(2\pi i\dot{z}_{ij})\;(i\leq j),\quad
\ddot{q}_{ij}:=\text{exp}(2\pi i\ddot{z}_{ij})\;(i<j)
\]
then 
\[f=\sum a_f(H)e^{2\pi i\text{tr}(HZ)}=\sum a_f(H)q^H
\]
may be considered as an element of the formal power series ring
\[
\mathbb{C}[\dot{q}_{ij}^{\pm 1},\ddot{q}_{ij}^{\pm 1}(i<j)][\![
\dot{q}_{11},\ldots ,\dot{q}_{nn}]\!].
\]
Let $R$ be a subring of $\mathbb{C}$. We define
\begin{align*}
&M_k(U_n(\mathcal{O}_{\boldsymbol{K}}),\nu)_R\\
&:=\{ f=\sum a_f(H)q^H\in M_k(U_n(\mathcal{O}_{\boldsymbol{K}}),\nu)\,|\, a_f(H)\in R\;
{\rm for}\;{\rm all}
\; H\in \Lambda_n(\boldsymbol{K})\}
\end{align*}
and
\[
M_k^{(s)}(U_n(\mathcal{O}_{\boldsymbol{K}}),\nu)_R:=M_k(U_n(\mathcal{O}_{\boldsymbol{K}}),
\nu)_R\cap 
M_k^{(s)}(U_n(\mathcal{O}_{\boldsymbol{K}}),\nu).
\]
So we may consider the inclusion:
\[
M_k(U_n(\mathcal{O}_{\boldsymbol{K}}),\nu)_R\subset
R[\dot{q}_{ij}^{\pm 1},\ddot{q}_{ij}^{\pm 1}(i<j)][\![
\dot{q}_{11},\ldots ,\dot{q}_{nn}]\!].
\]

\section{Siegel modular forms}
\label{sec:3}
In this section we introduce some results concerning Siegel modular forms which
are needed in later sections.
\subsection{Definition and notation}
\label{sec:3.1}
Let $M_k(\Gamma_n)$ denote the space of Siegel modular forms of weight $k$
$(\in\mathbb{Z})$ for the Siegel modular group $\Gamma_n:=Sp_n(\mathbb{Z})$
and $S_k(\Gamma_n)$ the subspace of cusp forms. Any Siegel modular form
$f(Z)$ in $M_k(\Gamma_n)$  has a Fourier expansion of the form
\[
f(Z)=\sum_{0\leq T\in\Lambda_n}a_f(T)e^{2\pi i\text{tr}(TZ)},
\]
where
\[
\Lambda_n=Sym_n^*(\mathbb{Z})
:=\{ T=(t_{ij})\in Sym_n(\mathbb{Q})\;|\; t_{ii},\;2t_{ij}\in\mathbb{Z}\; \}
\]
(the lattice in $Sym_n(\mathbb{R})$ of half-integral, symmetric matrices).\\
Taking $q_{ij}:=\text{exp}(2\pi iz_{ij})$ with $Z=(z_{ij})\in\mathbb{H}_n$, we write
\[
q^T:=e^{2\pi i\text{tr}(TZ)}=\prod_{1\leq i<j\leq n}
q_{ij}^{2t_{ij}}\prod_{i=1}^nq_{ii}^{t_{ii}}.
\]
Using this notation, we obtain \textit{the generalized
q-expansion}:
\begin{align*}
f=\sum_{0\leq T\in\Lambda_n}a_F(T)\,q^T&=\sum_{t_i}\left( \sum_{t_{ij}} 
a_f(T)\prod_{i<j}q_{ij}^{2t_{ij}}\right)\prod_{i=1}^nq_{ii}^{t_{ii}}\\
&\in\mathbb{C}[q_{ij}^{-1},q_{ij}][\![ q_{11},\ldots, q_{nn} ]\!].
\end{align*}
For any subring $R\subset\mathbb{C}$, we adopt the notation,
\begin{align*}
& M_k(\Gamma_n)_R:=\{ f=\sum_{T\in\Lambda_n}a_f(T)q^T\;|\;
a_f(T)\in R\;(\forall T\in\Lambda_n)\;\},\\
& S_k(\Gamma_n)_R:=M_k(\Gamma_n)_R\cap S_k(\Gamma_n).
\end{align*}
Any element $f\in M_k(\Gamma_n)_R$ can be regarded as an element of
\[
R[q_{ij}^{-1},q_{ij}][\![ q_{11},\ldots, q_{nn}]\!].
\]

\subsection{Siegel modular forms of degree 2}
\label{sec:3.2} 
For any Siegel modular form
\[
f=\sum a_f(T)q^T \in M_k(\Gamma_n)_{\mathbb{Z}_{(p)}},
\]
there exists a formal power series correspondence,
\[
\widetilde{f}:=\sum\widetilde{a_f(T)}q^T\in
\mathbb{F}_p[q_{ij}^{-1},q_{ij}][\![q_{11},\ldots, q_{nn}]\!].
\] 
where $\widetilde{a_f(T)}$ denotes the reduction modulo $p$ of $a_f(T)$. We define
\begin{align*}
\widetilde{M}_k(\Gamma_n)_p:&=\{\widetilde{f}=\sum\widetilde{a_f(T)}q^T|
f\in M_k(\Gamma_n)_{\mathbb{Z}_{(p)}}\}\\
&\subset \mathbb{F}_p[q_{ij}^{-1},q_{ij}][\![ q_{11},\ldots, q_{nn}]\!].
\end{align*}
\begin{Def}
The algebra
\[
\widetilde{M}(\Gamma_n)_p:=\sum_{k\in\mathbb{Z}} \widetilde{M}_k(\Gamma_n)_p
\quad (resp. \;\;\widetilde{M}^{(e)}(\Gamma_n)_p:=\sum_{k\in 2\mathbb{Z}} \widetilde{M}_k(\Gamma_n)_p)
\]
is called \textit{the algebra of Siegel modular forms mod p}
(resp. \textit{the algebra of Siegel modular forms mod p of even weight}).
\end{Def}
The structure of $\widetilde{M}(\Gamma_2)_p$ was studied by Nagaoka~\cite{Na}. Here we introduce the structure theorem of $\widetilde{M}^{(e)}(\Gamma_2)_p$ for the cases $p\geq 5$.

We take the usual generators $G_4$, $G_6$, $X_{10}$ and $X_{12}$ of the graded ring of Siegl modular forms of even weight for $\Gamma _2$. Here, $G_k$ ($k=4$, $6$) is normalized by $a_{G_k}(O_2)=1$ and $X_k$ ($k=10$, $12$) is normalized by $a_{X_k}\left(
\begin{pmatrix}1 & \tfrac{1}{2} \\ \tfrac{1}{2} & 1 \end{pmatrix}\right)
=1$. As in \cite{Na}, it is known that the graded ring over $\mathbb{Z}_{(p)}$ of Siegl modular forms of even weight with $p$-integral Fourier coefficients for $\Gamma _2$ is genrated by $G_4$, $G_6$, $X_{10}$ and $X_{12}$. 

\begin{Thm}[Nagaoka~\cite{Na}]
\label{mod p Siegel}
Assume that $p\geq 5$. There exists a Siegel modular form $F_{p-1}\in M_{p-1}(\Gamma_2)_{\mathbb{Z}_{(p)}}$ such that
\[
F_{p-1} \equiv 1 \bmod{p}
\]
and one has
\[
\widetilde{M}^{(e)}(\Gamma_2) \cong \mathbb{F}_p[x_1,x_2,x_3,x_4]/(\widetilde{A}-1), 
\]
where $(\widetilde{A}-1)$ is a principal ideal generated by $\widetilde{A}-1$ and $A\in\mathbb{Z}_{(p)}[x_1,x_2,x_3,x_4]$ is defined by
\[
F_{p-1}=A(G_4,G_6,X_{10},X_{12}).
\]
\end{Thm}

The following theorem is a generalization of Sturm's theorem in the case of Siegel modular forms of degree $2$. 
\begin{Thm}[Kikuta~\cite{Ki}]
\label{Ki}
Let $p$ be a prime with $p\ge 5$, $k$ an even integer and $\Gamma $ an arbitrary congruence subgroup of $\Gamma _2$. We put $t:=[\Gamma _2:\Gamma ]$. If $f\in M_k(\Gamma )_{\mathbb{Z}_{(p)}}$ satisfies that $a_f\left(\begin{pmatrix} m & \frac{r}{2} \\ \frac{r}{2} & n \end{pmatrix}\right)\equiv 0$ mod $p$ for all $m$, $n$, $r$ with $0\le m\le (kt+2)/10$, $0\le n\le (3kt+1)/30$ and $mn-r^2/4\ge 0$, then we have $f\equiv 0$ mod $p$. 
\end{Thm}
\begin{Rem}
(1) In general, $m$, $n$, $r$ are not integers. If $\Gamma $ is a congruence subgroup of level $N$, then the above $\begin{pmatrix} m & \frac{r}{2} \\ \frac{r}{2} & n \end{pmatrix}$ is an element of $\frac{1}{N}\Lambda _2$. \\
(2) In fact, we need only to check the case $n\le m$. 
\end{Rem}

Ichikawa \cite{Ich} showed that Serre's $p$-adic weight is well defined in the case of Siegel modular forms. We introduce the result in the general degree case. 
\begin{Thm}[Ichikawa \cite{Ich}]
\label{Ich}
Let $p$ be a prime with $p>n+3$ or $p\equiv 1$ mod $4$ and $f\in M_k(\Gamma _n)_{\mathbb{Z}_{(p)}}$, $g\in M_{k'}(\Gamma _n)_{\mathbb{Z}_{(p)}}$. If $f\equiv g$ mod $p^l$ and $f\not \equiv 0$ mod $p$, then we have $k\equiv k'$ mod $(p-1)p^{l-1}$. 
\end{Thm}

\section{Hermitian modular forms of degree 2}
\label{sec:4}

 In this section, we deal with Hermitian modular forms of degree 2. For more detail, we refer to \cite{Der,D-K}.

We consider the Hermitian Eisenstein series of degree 2
\[
E_k(Z):=\sum_{M=\left(\begin{smallmatrix} * & * \\ C & D \end{smallmatrix}\right)}({\det}M)^{\frac{k}{2}}{\det}(CZ+D)^{-k},\quad 
Z\in\mathbb{H}_2,
\]
where $k>4$ is even and $M=\Big(\begin{smallmatrix} * & * \\ C & D \end{smallmatrix}\Big)$ runs over a set of representatives of $\left\{\Big(\begin{smallmatrix} * & * \\ O_2 & * \end{smallmatrix}\Big)\right\}\backslash U_2(\mathcal{O}_{\boldsymbol{K}})$. 
Then we have
\[
E_k\in M_k^{(s)}(U_2(\mathcal{O}_{\boldsymbol{K}}),{\det}^{-\frac{k}{2}}).
\]
Moreover $E_4\in M_4^{(s)}(U_2(\mathcal{O}_{\boldsymbol{K}}),{\det}^{-2})$ is constructed by the Maass lift (\cite{Kr}).


In the rest of this paper, we mainly deal with $\boldsymbol{K}=\mathbb{Q}(\sqrt{-1})$ or $\mathbb{Q}(\sqrt{-3})$ and 
\begin{align*}
{\nu_k=}
\begin{cases}
{\det}^{k/2} & \text{for}\;\; \boldsymbol{K}=\mathbb{Q}(\sqrt{-1}),\\
{\det}^k & \text{for}\;\; \boldsymbol{K}=\mathbb{Q}(\sqrt{-3}).
\end{cases}
\end{align*}
We remark that $\nu _k$ is a trivial character if $\# \mathcal{O}_{\boldsymbol{K}}^{\times}|k$. This fact follows from that there exists a unit $\varepsilon \in \mathcal{O}_{\boldsymbol{K}}^{\times}$ such that $\det M =\varepsilon ^2$ if $M\in U_2(\mathcal{O}_{\boldsymbol{K}})$. The graded rings over $\mathbb{C}$ of these Hermitian modular forms are studied by Dern and Krieg~\cite{D-K} and \cite{Der}. 

\begin{Prop}[cf.~\cite{Der}]
\label{coeff}
Let $f\in M_k(U_2(\mathcal{O}_{\boldsymbol{K}}),\nu _k)$. There is the following relation between the Fourier coefficients:
\begin{align*}
a_{f| _{\mathbb{S}_2}}\left(\begin{pmatrix} m & \frac{r}{2} \\ \frac{r}{2} & n \end{pmatrix}\right)=\sum _{\substack{\alpha \in {\mathcal D}^{-1}_{\boldsymbol K} \\ \alpha +\overline{\alpha }=r \\ mn-N(\alpha )\ge 0}}a_f\left(\begin{pmatrix} m & \alpha \\ \overline{\alpha } & n \end{pmatrix}\right), 
\end{align*}
where ${\mathcal D}^{-1}_{\boldsymbol K}$ is the inverse different ${\mathcal D}^{-1}_{\boldsymbol K}:=\frac{{\mathcal O_{{\boldsymbol K}}}}{\sqrt{d_{{\boldsymbol K}}}}$. 
\end{Prop}

We define a lexicographical order for the different elements 
\begin{align*}
H=\begin{pmatrix} m & \frac{a+bi}{2} \\ \frac{a-bi}{2} & n \end{pmatrix},\quad H'=\begin{pmatrix} m' & \frac{a'+b'i}{2} \\ \frac{a'-b'i}{2} & n' \end{pmatrix} \in \Lambda_2({\boldsymbol K})
\end{align*} 
by
\begin{align*}
&H> H'\Longleftrightarrow (1)\ {\rm tr}(H)>{\rm tr}(H') \quad {\rm or}\quad  (2)\ {\rm tr}(H)={\rm tr}(H'),\ m>m' \\ 
&\quad {\rm or}\quad (3)\ {\rm tr}(H)={\rm tr}(H'),\ m=m',\ a> a' \\
&\quad {\rm or}\quad (4)\ {\rm tr}(H)={\rm tr}(H'),\ m=m',\ a=a',\ b>b'. 
\end{align*}
Let $p$ be a prime and $f\in M_k(U_2(\mathcal{O}_{\boldsymbol{K}}),\nu _k)_{\mathbb{Z}_{(p)}}$. We define an order of $f$ by 
\begin{align*}
{\rm ord}_p(f)=\min \{H|a_f(H) \not \equiv 0 \bmod{p}\}, 
\end{align*}
where the ``minimum'' is defined in the sense of the order defined above. If $f\equiv 0$ mod $p$, then we define ${\rm ord}_p(f)=\infty $. 

It is not difficult to see the following property. 
\begin{Lem}
\label{Lem}
One has 
\begin{align*}
{\rm ord}_p(fg)={\rm ord}_p(f)+{\rm ord}_p(g). 
\end{align*}
\end{Lem}

\begin{Prop}[Kikuta-Nagaoka~\cite{Ki-Na}]
\label{Ki-Na0}
Let ${\boldsymbol K}=\mathbb{Q}(\sqrt{-1})$ or $\mathbb{Q}(\sqrt{-3})$ and $p$ be a prime with $p\geq 5$. Then there exists a Hermitian modular form $F_{p-1}\in M_{p-1}^{(s)}(U_2({\mathcal O}_{\boldsymbol K}),\nu_{p-1})_{\mathbb{Z}_{(p)}}$ such that 
\[
F_{p-1}\equiv 1 \bmod{p}. 
\]
\end{Prop}

\begin{Thm}[Kikuta-Nagaoka~\cite{Ki-Na}]
\label{Ki-Na}
Let ${\boldsymbol K}=\mathbb{Q}(\sqrt{-1})$ and assume that $p\geq 5$. There exist cusp forms $\chi_{8}$, $F_{10}$, $F_{12}$ of respective weight $8$, $10$, $12$ satisfying the following. \\
\noindent
(1) $E_4$, $E_6$, $\chi_{8}$, $F_{10}$ and $F_{12}$ are algebraically independent and
\[
\bigoplus _{0\leq k\in 2\mathbb{Z}}M_{k}(U_2({\mathcal O}_{\boldsymbol K}),
\nu_k)^{sym }_{\mathbb{Z}_{(p)}}=\mathbb{Z}_{(p)}[E_4, E_6, \chi_8, F_{10}, F_{12}]. 
\]
In other words, if $f\in M_k^{(s)}(U_2(\mathcal{O}_{\boldsymbol K}),\nu_k)_{\mathbb{Z}_{(p)}}$ ($k$\;:\;even), then there exists a unique polynomial $P(x_1,x_2,x_3,x_4,x_5)\in\mathbb{Z}_{(p)}[x_1,x_2,x_3,x_4,x_5]$ such that 
\[
f=P(E_4, E_6, \chi_{8}, F_{10}, F_{12}). 
\]
\noindent 
(2) $E_4|_{\mathbb{S}_2}=G_4$,\quad $E_6|_{\mathbb{S}_2}=G_6$,\quad $\chi_8|_{\mathbb{S}_2} 
\equiv 0$,
\quad $F_{10}|_{\mathbb{S}_2}=6X_{10}$,\quad $F_{12}|_{\mathbb{S}_2}=X_{12}$. \\\noindent 
(3) $a_{\chi _8}\left(\begin{pmatrix} 1 & \frac{-1-i}{2} \\ \frac{-1+i}{2} & 1 \end{pmatrix}\right)=1$, namely ${\rm ord}_p(\chi _8)=\begin{pmatrix} 1 & \frac{-1-i}{2} \\ \frac{-1+i}{2} & 1 \end{pmatrix}$. 
\end{Thm}

\begin{Thm}[Kikuta-Nagaoka~\cite{Ki-Na}]
Let ${\boldsymbol K}=\mathbb{Q}(\sqrt{-3})$ and assume that $p\geq 5$. There exist cusp forms $F_{10}$, $F_{12}$, $\chi _{18}$ of respective weight $10$, $12$, $18$ satisfying the following. \\
\noindent 
(1) $E_4$, $E_6$, $F_{10}$, $F_{12}$, $\chi _{18}$ are algebraically independent and 
\[
\bigoplus _{0\leq k\in 2\mathbb{Z}}M_{k}(U_2({\mathcal O}_{\boldsymbol K}),
\nu_k)^{sym }_{\mathbb{Z}_{(p)}}=\mathbb{Z}_{(p)}[E_4, E_6, F_{10}, F_{12}, \chi _{18}]. 
\]
In other words, if $f\in M_k^{(s)}(U_2(\mathcal{O}_{\boldsymbol K}),\nu_k)_{\mathbb{Z}_{(p)}}$ ($k$\;:\;even), then there exists a unique polynomial $P(x_1,x_2,x_3,x_4,x_5)\in\mathbb{Z}_{(p)}[x_1,x_2,x_3,x_4,x_5]$ such that 
\[
f=P(E_4, E_6, F_{10}, F_{12}, \chi _{18}). 
\]

\noindent 
(2) $E_4|_{\mathbb{S}_2}=G_4$,\quad $E_6|_{\mathbb{S}_2}=G_6$,\quad $F_{10}|_{\mathbb{S}_2}=2X_{10}$,\quad $F_{12}|_{\mathbb{S}_2}=2X_{12}$,\quad $\chi _{18}|_{\mathbb{S}_2}\equiv 0$. \\
\noindent
(3) $a_{\chi _{18}}\left(\begin{pmatrix} 2 & \frac{2i}{\sqrt{3}} \\ -\frac{2i}{\sqrt{3}} & 2 \end{pmatrix}\right)=1$, namely ${\rm ord}_p(\chi _{18})=\begin{pmatrix} 2 & * \\ * & 2 \end{pmatrix}$ for some $*$. 
\end{Thm}
For explicit expressions of above generators, see \cite{Ki-Na}.


\section{Proofs}
\label{sec:5}
In this section, we shall prove our theorems. However, since the proof is similar, we prove only the case ${\boldsymbol K}=\mathbb{Q}(\sqrt{-1})$. 
\subsection{Proof of Theorem~\ref{ThmM}}
\label{subsec:5-1}
By Theorem~\ref{Ki-Na} (1), we can write $f$ in the form 
\begin{align*}
f=P(E_4,E_6,F_{10},F_{12})+\chi _8g, 
\end{align*}
where $P$ is a four variables polynomial over $\mathbb{Z}_{(p)}$ and $g\in M_{k-8}^{(s)}(U_2({\mathcal O}_{\boldsymbol K}),\nu _{k-8})_{\mathbb{Z}_{(p)}}$. Restricting both sides to $\mathbb{S}_2$, we obtain $f|_{\mathbb{S}_2}=P(G_4,G_6,6X_{10},X_{12})$ because of Theorem~\ref{Ki-Na} (2). Proposition~\ref{coeff} implies that $a_{f|_{\mathbb{S}_2}}\left(\begin{pmatrix} m & * \\ * & n \end{pmatrix}\right)\equiv 0$ mod $p$ for all $m$, $n$ with $0\le m \le (5k+8)/40$ and $0\le n\le (15k+4)/120$. Note that $(k+2)/10<(5k+8)/40$ and $(3k+1)/30<(15k+4)/120$. By Theorem~\ref{Ki}, we have $f|_{\mathbb{S}_2}=P(G_4,G_6,6X_{10},X_{12})\equiv 0$ mod $p$. By Theorem~\ref{mod p Siegel}, we have $\widetilde{P}(x_1,x_2,6x_3,x_4) \in (\widetilde{A}-1)$. Since $\widetilde{P}$ is an isobaric polynomial and $6\in \mathbb{Z}_{(p)}^{\times}$, we obtain $\widetilde{P}=0$. Accordingly $P(E_4,E_6,F_{10},F_{12})\equiv 0$ mod $p$. Hence $a_{\chi _8g}\left(\begin{pmatrix} m & * \\ * & n \end{pmatrix}\right)\equiv 0$ for all $m$, $n$ with $0\le m\le (5k+8)/40$, $0\le n\le (15k+4)/120$. 

Now we assume that the diagonal of ${\rm ord}_p(g\chi _8)$ are $m_0$ and $n_0$: 
\begin{align*}
{\rm ord}_p(g\chi _8)=\begin{pmatrix} m_0 & * \\ * & n_0 \end{pmatrix}. 
\end{align*}
Note that $m_0>(5k+8)/40$ and $n_0>(15k+4)/120$. By Lemma~\ref{Lem}, we have 
\begin{align*}
{\rm ord}_p(g)=\begin{pmatrix} m_0-1 & * \\ * & n_0-1 \end{pmatrix}. 
\end{align*}
Then we see that $m_0-1>(5k+8)/40-1=(5(k-8)+8)/40$ and $n_0-1> (15k+4)/120-1=(15(k-8)+4)/120$. This means that $a_{g}\left(\begin{pmatrix} m & * \\ * & n \end{pmatrix}\right)\equiv 0$ mod $p$ for all $m$, $n$ with $0\le m\le (5(k-8)+8)/40$, $0\le n\le (15(k-8)+4)/120$. Using an inductive argument on the weight, we see that $g\equiv 0$ mod $p$. This completes the proof of Theorem \ref{ThmM}. $\square$

\begin{Rem}
This proof depends on the congruence criterion for Siegel modular forms of degree $2$. Accordingly, if it's criterion is more sharp, then our results would be more sharp. 
\end{Rem}


\subsection{Proof of Theorem~\ref{ThmM3}}

By Theorem~\ref{Ki-Na} (1), we can write $f$ and $g$ in the forms 
\begin{align}
\label{eq1}
f=P_1(E_4,E_6,F_{10},F_{12})+\chi _8f_1,\quad g=Q_1(E_4,E_6,F_{10},F_{12})+\chi _8g_1
\end{align}
where $P_1$ and $Q_1$ are polynomials of four variables over $\mathbb{Z}_{(p)}$ and $f_1\in M_{k-8}^{(s)}(U_2({\mathcal O}_{\boldsymbol K}),\nu _{k-8})_{\mathbb{Z}_{(p)}}$ and $g_1\in M_{k'-8}^{(s)}(U_2({\mathcal O})_{\boldsymbol K},\nu _{k'-8})_{\mathbb{Z}_{(p)}}$. By restricting two equations in (\ref{eq1}) to $\mathbb{S}_2$ and Proposition~\ref{coeff}, we have
\begin{align*}
P_1(G_4,G_6,6X_{10},X_{12})\equiv Q_1(G_4,G_6,6X_{10},X_{12}) \bmod{p^l}.
\end{align*}
If $P_1(G_4,G_6,6X_{10},X_{12})\not \equiv 0$ mod $p$, then we obtain $k\equiv k'$ mod $(p-1)p^{l-1}$ by Theorem~\ref{Ich}. Hence we take the largest number $s\le l$ such that $P_1(G_4,G_6,6X_{10},X_{12})\equiv Q_1(G_4,G_6,6X_{10},X_{12})\equiv 0$ mod $p^s$. By Theorem~\ref{mod p Siegel}, we obtain $\tilde{P}_1(x_1,x_2,6x_3,x_4)\in (\tilde{A}-1)$. Since $\tilde{P}_1$ is an isobaric polynomial, we have $\tilde{P}_1=0$, namely $P_1\equiv 0$ mod $p$ as a polynomial. We obtain also $P_1\equiv 0$ mod $p^{s}$ inductively. We set
\begin{align*}
&P_1':=\frac{1}{p^s}P_1, \quad Q_1':=\frac{1}{p^s}Q_1\in \mathbb{Z}_{(p)}[x_1,x_2,x_3,x_4]. 
\end{align*}
By the property of $s$, we have
\begin{align*}
&P_1'(G_4,G_6,6X_{10},X_{12})\equiv Q_1'(G_4,G_6,6X_{10},X_{12}) \bmod{p^{l-s}},\\
 &P_1'(G_4,G_6,6X_{10},X_{12})\not \equiv 0 \bmod{p}. 
\end{align*}
By Theorem~\ref{Ich}, we obtain $k\equiv k'$ mod $(p-1)p^{l-s-1}$. Now we may assume that $k=k'+p^{l-s-1}(p-1)a$. Taking multiplication of $F_{p-1}$ which is given in Proposition~\ref{Ki-Na0}, we have 
\begin{align}
\label{eq2}
 Q_1'(E_4,E_6,F_{10},F_{12})\equiv Q_1'(E_4,E_6,F_{10},F_{12})F_{p-1}^{ap^{(l-s-1)}} \bmod{p^{l-s}}. 
\end{align}
By Theorem~\ref{Ki-Na} (1), we can write 
\begin{align*}
Q_1'(E_4,E_6,F_{10},F_{12})F_{p-1}^{ap^{(l-s-1)}}=R_1(E_4,E_6,F_{10},F_{12})+\chi _8 h_1, 
\end{align*}
where $R_1$ is a four variables polynomial over $\mathbb{Z}_{(p)}$ and $h\in M_{k-8}(U_2({\mathcal O})_{\boldsymbol K},\nu _{k-8})$. Note that $P_1'(E_4,E_6,F_{10},F_{12})$, $R_1(E_4,E_6,F_{10},F_{12})+\chi _8 h_1 \in M_k^{(s)}(U_2({\mathcal O}_{\boldsymbol K}),\nu _{k})_{\mathbb{Z}_{(p)}}$. Since (\ref{eq2}) and by using Theorem~\ref{Ki-Na} (2), we have 
\begin{align*}
P'_1(G_4,G_6,6X_{10},X_{12})\equiv R_1(G_4,G_6,6X_{10},X_{12}) \bmod{p^{l-s}}.
\end{align*}
Since the weights of both sides are same, by using Theorem~\ref{mod p Siegel} repeatedly, we get $P_1'(E_4,E_6,F_{10},F_{12})\equiv R_1(E_4,E_6,F_{10},F_{12}) \bmod{p^{l-s}}$. It follows that 
\begin{align*}
P_1(E_4,E_6,F_{10},F_{12})\equiv p^sR_1(E_4,E_6,F_{10},F_{12})\equiv Q_1(E_4,E_6,F_{10},F_{12})-p^s\chi _8 h_1 \bmod{p^{l}}. 
\end{align*}
By the assumption that $f\equiv g$ mod $p^l$, we have
\begin{align*}
\chi _8(f_1-h_1 p^s) \equiv \chi _8 g_1 \bmod{p^{l}}. 
\end{align*}
Now we substitute $f_1-h_1 p^s$ to $f_1$. Since $f\not \equiv 0$ mod $p$, it must be that $\chi _8f_1\not \equiv 0$ mod $p$. By Lemma~\ref{Lem}, we see also that $f_1\not \equiv 0$ mod $p$ and $f_1\equiv g_1$ mod $p^l$. We write again $f_1$ and $g_1$ in the forms 
\begin{align*}
f_1=P_2(E_4,E_6,F_{10},F_{12})+\chi _8f_2,\quad  g_1=Q_2(E_4,E_6,F_{10},F_{12})+\chi _8g_2,
\end{align*}
where $P_2$ and $Q_2$ are polynomials of four variables over $\mathbb{Z}_{(p)}$, $f_2\in M_{k-16}^{(s)}(U_2({\mathcal O}_{\boldsymbol K}),\nu _{k-16})_{\mathbb{Z}_{(p)}}$ and $g_2\in M_{k'-16}^{(s)}(U_2({\mathcal O}_{\boldsymbol K}),\nu _{k'-16})_{\mathbb{Z}_{(p)}}$. By Theorem~\ref{Ich}, if $P_2(G_4,G_6,6X_{10},X_{12})\not \equiv 0$ mod $p$, then we obtain $k-8\equiv k'-8$ mod $(p-1)p^{l-1}$, namely $k\equiv k'$ mod $(p-1)p^{l-1}$. If otherwise, by continuing the above argument repeatedly, we find a natural number $t$ such that $g_t|_{\mathbb{S}_2}\not \equiv 0$ mod $p$. In fact, it suffices to continue until $t$ satisfying $0\le k'-8t<8$. Applying Theorem~\ref{Ich} to $f_t$ and $g_t$, we obtain that $k-8t\equiv k'-8t$ mod $(p-1)p^{l-1}$, namely $k\equiv k'$ mod $(p-1)p^{l-1}$. This completes the proof of Theorem~\ref{ThmM3}. $\square$



\noindent
T. Kikuta \\
Department of Mathematics Kinki University Higashi-Osaka, 577-8502 Osaka, Japan \\
Tel.: +81-6-6721-2332\\
Fax: +81-6-6727-4301\\
E-mail: kikuta84@gmail.com
\end{document}